\documentclass[preprint,3p,12pt,pdf]{elsarticle}

\usepackage{lineno,hyperref}
\modulolinenumbers[5]

\usepackage{hyperref}

\usepackage{epstopdf}

\usepackage{tikz}
\usetikzlibrary{arrows}

\usepackage{pst-node}
\usepackage{tikz-cd}

\usepackage[shellescape]{gmp}

\usepackage{mathrsfs}
\usepackage{amssymb}
\usepackage{amsfonts}
\usepackage{latexsym}
\usepackage{amsmath}
\usepackage{xcolor}
\usepackage{wrapfig}
\usepackage{floatflt}

\usepackage{mathtools}
\usepackage{extarrows}

\usepackage{graphicx}
\def\lb{\label}

\newcommand{\er}[1]{\textrm{(\ref{#1})}}

\newtheorem{theorem}{\bf Theorem}[section]

\newtheorem{definition}[theorem]{\bf Definition}

  \def\cA{{\mathcal A}}       \def\mA{{\mathscr A}}
   \def\cB{{\mathcal B}}       \def\mB{{\mathscr B}}
  \def\cC{{\mathcal C}}       \def\mC{{\mathscr C}}
  \def\cD{{\mathcal D}}       
\def\d{\delta}         
         \def\mF{{\mathscr F}}
           
          \def\mH{{\mathscr H}}
    \def\cI{{\mathcal I}}       
           
 \def\cK{{\mathcal K}}       \def\mK{{\mathscr K}}
         
 \def\cM{{\mathcal M}}       \def\mM{{\mathscr M}}

\def\s{\sigma}         
  \def\cS{{\mathcal S}}       
           
    \def\cU{{\mathcal U}}       \def\mU{{\mathscr U}}

       \def\vp{\varphi}    

\def\Z{{\mathbb Z}}    \def\R{{\mathbb R}}   \def\C{{\mathbb C}}    
\def\T{{\mathbb T}}    \def\N{{\mathbb N}}   

\def\lt{\biggl}                  \def\rt{\biggr}
\def\ol{\overline}               \def\wt{\widetilde}
\def\no{\noindent}

\let\ge\geqslant                 \let\le\leqslant

\def\iy{\infty}
               
\def\ss{\subset}                 \def\ts{\times}
                \def\os{\oplus}
                 \def\ev{\equiv}
        
\def\el2{\ell^{\,2}}             \def\1{1\!\!1}

\def\dim{\mathop{\mathrm{dim}}\nolimits}

\def\BBox{\hspace{1mm}\vrule height6pt width5.5pt depth0pt \hspace{6pt}}

\let\ge\geqslant
\let\le\leqslant

\newcommand{\ca}{\begin{cases}}
\newcommand{\ac}{\end{cases}}
\newcommand{\ma}{\begin{pmatrix}}
\newcommand{\am}{\end{pmatrix}}
\def\eq{\begin{equation}}
\def\qe{\end{equation}}
\def\[{\begin{equation}}
\def\]{\end{equation}}

\def\BBox{\hspace{1mm}\vrule height6pt width5.5pt depth0pt \hspace{6pt}}

\begin{document}

\begin{frontmatter}



\title{Classification of integro-differential $C^*$-algebras}

\date{\today}

\author
{Anton A. Kutsenko}

\address{Jacobs University, 28759 Bremen, Germany; email: akucenko@gmail.com}

\begin{abstract}
The integro-differential algebra $\mathscr{F}_{N,M}$ is the $C^*$-algebra generated by the following operators acting on $L^2([0,1)^N\to\mathbb{C}^M)$: 1) operators of multiplication by bounded matrix-valued functions, 2) finite differential operators, 3) integral operators. Most of the operators and their approximations studying in physics belong to these algebras. We give a complete characterization of $\mathscr{F}_{N,M}$ in terms of its Bratteli diagram. In particular, we show that $\mathscr{F}_{N,M}$ does not depend on $M$ but depends on $N$. At the same time, it is known that differential algebras $\mathscr{H}_{N,M}$, generated by the operators 1) and 2), do not depend on both dimensions $N$ and $M$, they are all $*$-isomorphic to the universal UHF-algebra. We explicitly compute the Glimm-Bratteli symbols (for $\mH_{N,M}$ it was already computed earlier)  
$$
 \mathfrak{n}(\mathscr{F}_{N,M})=\prod_{n=1}^{\infty}\begin{pmatrix} n & 0 \\ n-1 & 1 \end{pmatrix}^{\otimes N}\begin{pmatrix}1 \\ 1 \end{pmatrix}^{\otimes N},\ \ \ \  \mathfrak{n}(\mathscr{H}_{N,M})=\prod_{n=1}^{\infty}n,
$$
which characterize completely the corresponding AF-algebras.  
\end{abstract}



\end{frontmatter}


{\section{Introduction}\lb{sec1}}

Discrete and continuous analogues of integro-differential algebras are actively used in various applications, for example, in the development of computer algorithms for symbolic and numerical solving of integro-differential equations, see, e.g., \cite{Ro1,GK01,Bu1,GRR}. On the other hand, differential algebras are closely related to the rotation $C^*$-algebras well studied in, e.g., \cite{PV1980,R1981,Y1986,EN2015}. In contrast to the rotation algebras, the integro-differential algebras contain operators of multiplication by discontinuous functions and integral operators. Nevertheless, the integro-differential algebras are AF-algebras and, hence, they admit a classification in terms of, e.g., the Bratteli diagrams.

{\section{Characterization of AF-algebras. Preliminary results.}\lb{sec2}}

Let us recall some facts about Bratteli diagrams. It is well known that any finite-dimensional $C^*$-algebra is $*$-isomorphic the direct sum of simple matrix algebras. Up to the order of terms, this direct sum is determined uniquely. It is convenient to use the following notation for finite-dimensional $C^*$-algebras. Let ${\bf p}=(p_j)_{j=1}^n\in\N^n$, then
\[\lb{prod}
 \mM({\bf p})=\C^{p_1\ts p_1}\os...\os\C^{p_n\ts p_n}.
\]
Any $*$-homomorphism from $\mM({\bf p})$ to $\mM({\bf q})$ with ${\bf p}\in\N^n$, ${\bf q}\in\N^m$ is internally (inside each $\C^{q_j\ts q_j}$) unitary  equivalent to some canonical $*$-homomorphism. Any canonical $*$-homomorphism is completely and uniquely determined by the matrix of multiplicities of partial embeddings ${\bf E}\in\Z_+^{m\ts n}$ (E-matrix) satisfying ${\bf E}(p_j)_{j=1}^n=(\wt q_j)_{j=1}^m$, where $\wt q_j\le q_j$. For example, the canonical $*$-homomorphism
$$
 \vp: \mM(2,2,3)\to \mM(4,4),\ \ {\bf A}\os{\bf B}\os{\bf C}\xmapsto[]{\ \vp\ }\ma {\bf A} & 0 \\ 0 & {\bf A} \am\os\ma  0 \am 
$$
has the E-matrix
$$
 {\bf E}=\ma 2 & 0 & 0 \\ 0 & 0 & 0 \am.
$$
For simplicity, we can write
$$
 \mM(2,2,3)\xrightarrow[]{\ {\bf E}\ }\mM(4,4).
$$
If a canonical $*$-homomorphism is unital then there are no zero rows in E-matrix, and we should replace the above mentioned condition $\wt q_j\le q_j$ with $\wt q_j=q_j$. For example, the unital  embedding
$$
 \mM(2,2,3)\xrightarrow[]{\ {\bf E}\ }\mM(4,5),\ \ {\bf E}=\ma 2 & 0 & 0 \\ 0 & 1 & 1 \am
$$
has the form
$$
 {\bf A}\os{\bf B}\os{\bf C}\xmapsto[]{}\ma {\bf A} & 0 \\ 0 & {\bf A} \am\os\ma {\bf B} & 0 \\ 0 & {\bf C} \am.
$$

The AF-algebra is a separable $C^*$-algebra, any finite subset of which can be approximated by a finite-dimensional $C^*$-sub-algebra. For convenience, we will consider unital AF-algebras only. This is not a restriction because the unitalization of an AF-algebra is obviously an AF-algebra. It is well known, that for any unital AF-algebra $\mA$ there is a family of nested finite-dimensional $C^*$-subalgebras $\mA_n\subseteq\mA$, satisfying
\[\lb{nest1}
 \C^{1\ts1}\cong\mA_0\subseteq\mA_1\subseteq\mA_2\subseteq...,\ \ \ \mA=\ol{\bigcup_{n=0}^{\infty}\mA_n}.
\] 
Since $\mA_n$ are nested finite-dimensional $C^*$-algebras, they are isomorphic to some canonical algebras $\mM({\bf p}_n)$, where ${\bf p}_n\in\N^{M_n}$, $M_n\in\N$, and the inclusions \er{nest1} can be written as
\[\lb{nest2}
 \mM({\bf p}_0)\xrightarrow[]{\ {\bf E}_0\ }\mM({\bf p}_1)\xrightarrow[]{\ {\bf E}_1\ }\mM({\bf p}_2)\xrightarrow[]{\ {\bf E}_2\ }...,
\]
where ${\bf p}_0=1$, $M_0=1$, and ${\bf E}_n\in\Z_+^{M_{n+1}\ts M_n}$. Moreover, due to the unital embeddings $\mM({\bf p}_n)\subseteq\mM({\bf p}_{n+1})$, all E-matrices have no zero rows and columns and they satisfy ${\bf E}_n{\bf p}_n={\bf p}_{n+1}$. Because ${\bf p}_0=1$, we obtain
\[\lb{nest3}
{\bf p}_{n+1}={\bf E}_n...{\bf E}_1{\bf E}_0=\prod_{i=0}^n{\bf E}_n.
\]
We will always assume the right-to-left order in the product $\prod$. Using \er{nest2}-\er{nest3}, we conclude that the matrices $\{{\bf E}_n\}_{n=0}^{+\iy}$ determine completely the structure of the unital AF-algebra $\mA$. It is useful to note that the choice of E-matrices is not unique. For example, E-matrices $\{{\bf E}_n'\}_{n=0}^{+\iy}$, where ${\bf E}_n'={\bf E}_{2n+1}{\bf E}_{2n}$, determine the same algebra $\mA$. This is because the composition of two embeddings has the E-matrix equivalent to the product of E-matrices corresponding to the embeddings. It is possible to describe the class of all E-matrices determining the same unital AF-algebra.

\begin{definition}\lb{D1} Let $\mathfrak{E}$ be the set of sequences of matrices $\{{\bf E}_n\}_{n=0}^{\iy}$, where ${\bf E}_n\in\Z_+^{M_{n+1}\ts M_n}$ have no zero rows and columns, and $M_0=1$, $M_n\in\N$  are some positive integer numbers. Let us define the equivalence relation on $\mathfrak{E}$. Two sequences $\{{\bf A}_n\}_{n=0}^{\iy}\sim\{{\bf B}_n\}_{n=0}^{\iy}$ are equivalent if there is $\{{\bf C}_n\}_{n=0}^{\iy}\in\mathfrak{E}$ such that
\[\lb{D11}
 {\bf C}_0=\prod_{i=0}^{r_1-1}{\bf A}_i,\ \ {\bf C}_{2n-1}{\bf C}_{2n-2}=\prod_{i=m_{n-1}}^{m_{n}-1}{\bf B}_i,\ \ {\bf C}_{2n}{\bf C}_{2n-1}=\prod_{i=r_n}^{r_{n+1}-1}{\bf A}_i,\ \ n\ge1,
\]
where $0=r_0<r_1<r_2<...$ and $0=m_0<m_1<m_2<...$ are some monotonic sequences of integer numbers. The corresponding set of equivalence classes is denoted by $\mathfrak{E}:=\mathfrak{E}/\sim$.
\end{definition}

It is convenient to denote the equivalence classes as
$$
 \ol{\{{\bf E}_n\}_{n=0}^{\iy}}=\prod_{n=0}^{\infty}{\bf E}_n
$$
because, see \er{D11},
$$
 \prod_{n=0}^{\iy}{\bf A}_n=\prod_{n=0}^{\iy}\prod_{i=r_n}^{r_{n+1}-1}{\bf A}_i=\prod_{n=0}^{\iy}{\bf C}_n=\prod_{n=0}^{\iy}\prod_{i=m_n}^{m_{n+1}-1}{\bf B}_i=\prod_{n=0}^{\iy}{\bf B}_n.
$$
In other words, we can perform the standard manipulations in the product of matrices without leaving the equivalence class. Of course, the manipulations should not go beyond $\mathfrak{E}$, i.e. all the resulting matrices should have non-negative integer entries and should not have zero rows and zero columns.

Let $\mA$ be some unital AF-algebra. Following \er{nest1}-\er{nest3}, there is a Bratteli diagram $\{{\bf E}_n\}_{n=0}^{\iy}$ represented $\mA$. Let us define the mapping
\[\lb{GBmap}
 \mathfrak{n}:\mA\mapsto\prod_{n=0}^{\iy}{\bf E}_n.
\]
Because the Bratteli diagram is not unique, the correctness of the mapping $\mathfrak{n}$ should be checked. It is already done in the main structure theorem for Bratteli diagrams.

\begin{theorem}\lb{Tstruc}
i) The relation $\sim$ defined in \er{D11} is the equivalence relation. ii) Let $\mathfrak{A}$ be the set of classes of non-isomorphic unital AF-algebras. Then $\mathfrak{n}:\mathfrak{A}\to\mathfrak{E}$ is 1-1 mapping.
\end{theorem}

Note that the inverse mapping $\mathfrak{n}^{-1}$ has a more explicit form than $\mathfrak{n}$. For example, $\mathfrak{n}^{-1}(\prod_{n=0}^{\iy}{\bf E}_n)$ is the $C^*$-algebra $\mA$ given by the inductive limit \er{nest2}.

The proof of Theorem \ref{Tstruc} follows from the similar results formulated for the graphical representations of Bratteli diagrams, see, e.g., \cite{B1972,D1997}, and Theorem 3.4.4 in \cite{N2016}. The equivalence relation $\sim$ defined in \er{D11} is the analogue of telescopic transformations of Bratteli diagrams.
While the Bratteli diagram is not unique, it provides a kind of classification tool. Other types of classification of AF algebras, including the efficient K-theoretic Elliott classification, are discussed in \cite{E1976,EHS1980,RLL2000,AEG2015}. The infinite product $\mathfrak{n}(\mA)$ representing the Bratteli diagram for AF-algebra $\mA$ can be called as the Glimm-Bratteli symbol. Using supernatural symbols (numbers), when ${\bf E}_n$ are natural numbers in \er{GBmap}, J. Glimm provides the classification of uniformly hyper-finite algebras in \cite{G1960}.

Let us consider some examples of AF-algebras.

{\bf 1. Compact operators.} Let $\mK$ be the $C^*$-algebra of compact operators acting on a separable Hilbert space. Let $\mK_1={\rm Alg}(\mK,1)$ be its unitalization. It is well known that the Bratteli diagram for $\mK_1$ is

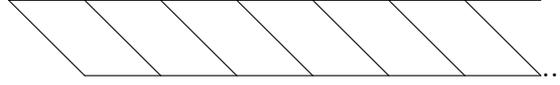
\qquad\qquad\qquad\qquad\qquad \begin{mpost*}
beginfig(1);
draw (3cm,1cm) -- (4cm,1cm);
draw (3cm,1cm) -- (4cm,0cm);
for i = 4cm step 1cm until 9cm:
  draw (i,0cm) -- (i+1cm,0cm); 
  draw (i,1cm) -- (i+1cm,0cm);
  draw (i,1cm) -- (i+1cm,1cm);  
endfor;
endfig;          
\end{mpost*}..,

\no where the nodes represent simple matrix sub-algebras, and the edges show the multiplicity of embedding: one line means the multiplicity equal to $1$. The first node is always $\C^{1\ts1}$. The dimensions of nodes are determined by the dimensions of nodes connected on the left and by the multiplicities of embedding. The corresponding Glimm-Bratteli symbol is
$$
 \mathfrak{n}(\mK_1)=\ma 1 & 0 \\ 1 & 1 \am^{\iy}\ma 1 \\ 1 \am.
$$
Combining terms in the infinite product, we can write the another form of the Glimm-Bratteli symbol
$$
 \mathfrak{n}(\mK_1)=\ma 1 & 0 \\ 1 & 1 \am^{\iy}\ma 1 \\ 1 \am=\lt(\prod_{n=1}^{\iy}\ma 1 & 0 \\ 1 & 1\am^n\rt)\ma 1 \\ 1 \am=\lt(\prod_{n=1}^{\iy}\ma 1 & 0 \\ n & 1\am\rt)\ma 1 \\ 1 \am,
$$ 
which leads to the labeled Bratteli diagram

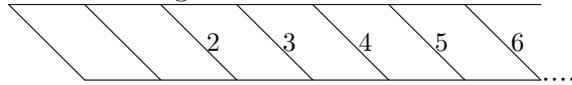
\qquad\qquad\qquad\qquad\qquad \begin{mpost*}
beginfig(2);
draw (3cm,1cm) -- (4cm,1cm);
draw (3cm,1cm) -- (4cm,0cm);
for i = 4cm step 1cm until 9cm:
  draw (i,0cm) -- (i+1cm,0cm); 
  draw (i,1cm) -- (i+1cm,0cm);
  draw (i,1cm) -- (i+1cm,1cm);  
endfor;
label.rt(btex $2$ etex, (5cm+0.5cm,0.5cm));
label.rt(btex $3$ etex, (6cm+0.5cm,0.5cm));
label.rt(btex $4$ etex, (7cm+0.5cm,0.5cm));
label.rt(btex $5$ etex, (8cm+0.5cm,0.5cm));
label.rt(btex $6$ etex, (9cm+0.5cm,0.5cm));
endfig;          
\end{mpost*}....

\no In the labeled Bratteli diagram, the edge numbers are the multiplicities of embedding. The multiplicity $1$ is usually omitted.

{\bf 2. CAR algebra.} For the CAR algebra $\mC$, which is a UHF-algebra, we have the Glimm-Bratteli symbol $\mathfrak{n}(\mC)=2^{\iy}$. At the same time,
$$
 \mathfrak{n}(\mC)=2^{\iy}=\lt(\ma 1 & 1\am\ma 1 \\ 1 \am\rt)^{\iy}=\lt(\ma 1 \\ 1 \am\ma 1 & 1\am\rt)^{\iy}\ma 1 \\ 1 \am=\ma 1 & 1  \\ 1 & 1\am^{\iy}\ma 1 \\ 1 \am.
$$
The corresponding Bratteli diagrams are as follows

\begin{mpost*}
beginfig(3);
draw (0cm,0cm) -- (3cm,0cm);
draw (0cm,0.1cm) -- (3cm,0.1cm);
endfig;
\end{mpost*}... \ \ \ \ $=$\ \ \ \ 
\begin{mpost*}
beginfig(4);
draw (0cm,0cm) -- (3cm,0cm);
label.top(btex $2$ etex, (0.5cm,0cm));
label.top(btex $2$ etex, (1.5cm,0cm));
label.top(btex $2$ etex, (2.5cm,0cm));
endfig;
\end{mpost*}... \ \ \ \ $=$\ \ \ \ 
\begin{mpost*}
beginfig(5);
draw (0cm,1cm) -- (1cm,1cm);
draw (0cm,1cm) -- (1cm,0cm);
for i = 1cm step 1cm until 3cm:
  draw (i,0cm) -- (i+1cm,0cm); 
  draw (i,1cm) -- (i+1cm,1cm);
  draw (i,1cm) -- (i+1cm,0cm); 
  draw (i,0cm) -- (i+1cm,1cm); 
endfor;
endfig;          
\end{mpost*}....

{\bf 3. Direct sum of AF algebras.} 
Above, we already used the notation $\os$ for the direct sum of matrix algebras and for their elements. We will use the same symbol in a little bit different context, namely for the direct sum of not necessarily square E-matrices. Suppose that $\mC=\mA\os\mB$ is the standard direct sum of two AF-algebras. 
If $\mathfrak{n}(\mA)=\prod_{n=0}^{\iy}{\bf A}_n$ and $\mathfrak{n}(\mB)=\prod_{n=0}^{\iy}{\bf B}_n$ then it can be shown that 
$$
 \mathfrak{n}(\mC)=\lt(\prod_{n=0}^{\iy}{\bf C}_n\rt)\ma 1 \\ 1 \am,\ \ where\ \ {\bf C}_n={\bf A}_n\os{\bf B}_n=\ma {\bf A}_n & {\bf 0} \\ {\bf 0} & {\bf B}_n \am.
$$

{\bf 4. Tensor product of AF-algebras.} It is useful to note the following property of the tensor product
$$
 \mM({\bf p})\otimes \mM({\bf q})=\mM({\bf p}\otimes{\bf q}),
$$
where
$$
 {\bf p}=(p_i)_{i=1}^n\in\N^n,\ \ \ {\bf q}=(q_j)_{j=1}^m\in\N^m,\ \ \ {\bf p}\otimes{\bf q}=(p_iq_j)_{i,j=1}^{n,m}\in\N^{nm}.
$$
Moreover, it is easy to check that if
$$
 \mM({\bf p}_1)\xrightarrow[]{\ {\bf E}_1\ }\mM({\bf q}_1),\ \ \ \mM({\bf p}_2)\xrightarrow[]{\ {\bf E}_2\ }\mM({\bf q}_2)
$$
then
$$
 \mM({\bf p}_1\otimes{\bf p}_2)\xrightarrow[]{\ {\bf E}_1\otimes{\bf E}_2\ }\mM({\bf q}_1\otimes{\bf q}_2),
$$
where the tensor product of matrices is defined in the standard way
$$
 (A_{i,j})\otimes(B_{r,s})=(C_{(i,r),(j,s)}),\ \ C_{(i,r),(j,s)}=A_{i,j}B_{r,s}.
$$
Hence, the standard tensor product $\mC=\mA\otimes\mB$ of two AF-algebras is AF-algebra which satisfies
$$
 \mathfrak{n}(\mC)=\prod_{n=0}^{\iy}{\bf A}_n\otimes{\bf B}_n,
$$
where $\mathfrak{n}(\mA)=\prod_{n=0}^{\iy}{\bf A}_n$, $\mathfrak{n}(\mB)=\prod_{n=0}^{\iy}{\bf B}_n$, and the tensor product of matrices is then
$$
 {\bf A}\otimes{\bf B}=\ma b_{11}{\bf A} & ... & b_{1N}{\bf A} \\ ... & ... & ... \\ b_{M1}{\bf A} & ... & b_{MH}{\bf A} \am,\ \ {\bf B}=\ma b_{11} & ... & b_{1N} \\ ... & ... & ... \\ b_{M1} & ... & b_{MN} \am.
$$
We will also use the following result. 

\begin{theorem}\lb{T2} Let $\{{\bf A}_n\}_{n=1}^{\iy}$ be a commutative (multiplicative) semigroup of square matrices with non-negative integer entries and with non-zero determinants. Let ${\bf A}_0$ be a matrix-column with positive integer entries such that ${\bf A}_1{\bf A}_0$ is defined. Let $\{{\bf B}_n\}_{n=1}^{\iy}\ss\{{\bf A}_n\}_{n=1}^{\iy}$ be a subset consisting of not necessarily different matrices satisfying the condition ($\pmb{\s}$): for any $p\in\N$ there are $r,s\in\N$ such that ${\bf A}_p{\bf A}_r=\prod_{i=1}^s{\bf B}_i$. Then
\[\lb{eqgroup}
 \mathfrak{n}^{-1}((\prod_{n=1}^{\iy}{\bf B}_n){\bf A}_0)\cong \mathfrak{n}^{-1}((\prod_{n=1}^{\iy}{\bf A}_{n}){\bf A}_0).
\]
Even if ($\pmb{\s}$) is not fulfilled, LHS in \er{eqgroup} is a sub-algebra of RHS.
\end{theorem}

%
%
{\bf Remark.} The universal UHF-algebra $\mU$ is the AF-algebra generated by the multiplicative semigroup of natural numbers
$$
 \mU=\mathfrak{n}^{-1}(\prod_{n=1}^{+\iy}n)=\mathfrak{n}^{-1}(\prod_{n=1}^{+\iy}p_n^{\iy})=\mathfrak{n}^{-1}(\prod_{n=1}^{+\iy}(p_1...p_n)^n)=\mathfrak{n}^{-1}(\prod_{n=1}^{+\iy}(p_1...p_n)),
$$ 
where $p_1=2$, $p_2=3$, $p_3=5$, ... are the prime numbers. Any UHF-algebra is a sub-algebra of $\mU$. The CAR-algebra is the UHF-algebra generated by the multiplicative semigroups $\{2^n:\ n\in m\N\}$ for any $m\in\N$.

There is another useful proposition describing non-isomorphic classes of AF-algebras.

\begin{theorem}\lb{Tiso}
Let $N,M\in\N$. Let $\{{\bf A}_n\}\ss\Z_+^{N\ts N}$, $\{{\bf B}_n\}\ss\Z_+^{M\ts M}$ be two sequences of matrices having non-zero determinants. Let ${\bf A}_0\in\Z_+^{N\ts 1}$, ${\bf B}_0\in\Z_+^{M\ts 1}$ be two matrix-columns without zero entries. If $N\ne M$ then
$\mathfrak{n}^{-1}(\prod_{n=0}^{\iy}{\bf A}_n)\not\cong\mathfrak{n}^{-1}(\prod_{n=0}^{\iy}{\bf B}_n)$ are non-isomorphic $C^*$-algebras.
\end{theorem}

Suppose that we have two AF-algebras $\mA_1$ and $\mA_2$ acting on Hilbert spaces $H_1$ and $H_2$ respectively. Suppose that $\mA_1\cong\mA_2$ are isomorphic to each other. When can we construct the unitary $\cU:H_1\to H_2$ inducing the $C^*$-algebra isomorphism, i.e. $\mA_1=\cU^{-1}\mA_2\cU$?

Let us start with the simple situation $\mA_1\cong\C^{n\ts n}$. Recall that we always consider unital algebras in this paper.  Then, there is a decomposition $H_1=\bigoplus_{i=1}^nH_{1i}$, where all $H_{1i}$ have the same dimension. We may think that all $H_{1i}=H_{11}$ are the same Hilbert space. Any operator $\mA_1\ni\cA\simeq(a_{ij})_{i,j=1}^N$ has the form $\mA_1=(a_{ij}\cI)_{i,j=1}^N:\bigoplus_{j=1}^NH_{11}\to\bigoplus_{i=1}^NH_{11}$, where $\cI:H_{11}\to H_{11}$ is the identity operator. We call any decomposition of $H$ satisfying the explained property as a decomposition associated with the simple algebra $\mA_1$. 

Taking a decomposition $H_2=\bigoplus_{i=1}^NH_{22}$ associated with $\mA_2$, we may state: {\it the isomorphism between $C^*$-algebras $\mA_1\cong\mA_2(\cong\C^{n\ts n})$ can be induced by a unitary $\cU:H_1\to H_2$ if and only if $\dim H_{11}=\dim H_{22}$}.

Suppose that we have a finite dimensional algebra
\[\lb{terms1}
 \mM({\bf p})=\C^{p_1\ts p_1}\os...\os\C^{p_N\ts p_N},\ \ {\bf p}=(p_i)_{i=1}^N
\]
acting on a separable Hilbert space $H$. Then there is a decomposition $H=\bigoplus_{i=1}^NH_i$ such that each direct term in \er{terms1} acts on the corresponding $H_i$. Hence, we can take a decomposition $H_i=\bigoplus_{j=1}^{p_i} H_{i1}$ associated with the corresponding simple direct term. Let us define ${\bf q}=(q_i)_{i=1}^N$, where $q_i=\dim H_{i1}$. To show the internal structure of the algebra $\mM({\bf p})$, we will write $\mM({\bf p};{\bf q})$. Let $\mM({\bf p}_1;{\bf q}_1)$ be some (unital) sub-algebra with the corresponding embedding
$$
 \mM({\bf p}_1;{\bf q}_1)\xrightarrow[]{\ {\bf E}\ }\mM({\bf p};{\bf q})\ \ \ ({\bf p}={\bf E}{\bf p}_1).
$$
Then, it is not difficult to check that the dimensions satisfy
$$
 {\bf q}^{\top}{\bf E}={\bf q}_1^{\top}.
$$
Since the dimensions can be infinite, we should specify the rules:
$$
 a+b=\iy\ if\ a=\iy\ or\ b=\iy,\ \ \iy\cdot0=0\cdot\iy=0.
$$
Any (unital) AF-algebra acting on a separable Hilbert space can be represented through its Bratteli diagram
\[\lb{nest2d}
\mM(1;q_0)\xrightarrow[]{\ {\bf E}_0\ }\mM({\bf p}_1;{\bf q}_1)\xrightarrow[]{\ {\bf E}_1\ }\mM({\bf p}_2;{\bf q}_2)\xrightarrow[]{\ {\bf E}_2\ }...,
\]
where the dimensions satisfy
\[\lb{d3}
 {\bf p}_n={\bf E}_{n-1}...{\bf E}_0,\ \ {\bf q}_{n-1}^{\top}={\bf q}_n^{\top}{\bf E}_{n-1},\ \ n\ge1.
\]
In order to include the dimensions ${\bf q}_n$, let us extend Definition \ref{D1}.
\begin{definition}\lb{D2} Let $\mathfrak{F}$ be the set of sequences of matrices and dimensions $\{{\bf E}_n,{\bf q}_n\}_{n=0}^{\iy}$, where ${\bf E}_n\in\Z_+^{M_{n+1}\ts M_n}$ have no zero rows and columns, and $M_0=1$, $M_n\in\N$  are some positive integer numbers. 
The dimensions ${\bf q}_n\in(\N\cup\iy)^{M_n}$ satisfy
\[\lb{cdim}
 {\bf q}_{n-1}^{\top}={\bf q}_n^{\top}{\bf E}_{n-1},\ \ n\ge1.
\]
Let us define the equivalence relation on $\mathfrak{F}$. Two sequences $\{{\bf A}_n,{\bf a}_n\}_{n=0}^{\iy}\sim\{{\bf B}_n,{\bf b}_n\}_{n=0}^{\iy}$ are equivalent if there is $\{{\bf C}_n,{\bf c}_n\}_{n=0}^{\iy}\in\mathfrak{F}$ such that
\[\lb{D21}
 {\bf C}_0=\prod_{i=0}^{r_1-1}{\bf A}_i,\ \ {\bf C}_{2n-1}{\bf C}_{2n-2}=\prod_{i=m_{n-1}}^{m_{n}-1}{\bf B}_i,\ \ {\bf C}_{2n}{\bf C}_{2n-1}=\prod_{i=r_n}^{r_{n+1}-1}{\bf A}_i,\ \ n\ge1
\]
and
\[\lb{D22}
 {\bf c}_0={\bf a}_0={\bf b}_0,\ \ {\bf c}_{2n-1}={\bf a}_{r_n},\ \ {\bf c}_{2n}={\bf b}_{m_n},\ \ n\ge1,
\]
where $0=r_0<r_1<r_2<...$ and $0=m_0<m_1<m_2<...$ are some monotonic sequences of integer numbers. The corresponding set of equivalence classes is denoted by $\mathfrak{F}:=\mathfrak{F}/\sim$.
\end{definition}

Now, we can partially complement the results of Theorem \ref{Tstruc}.

\begin{theorem}\lb{Fiso} Let $H$ be a separable Hilbert space. We call two $C^*$-algebras $\mA_1$ and $\mA_2$ unitary equivalent iff there is a unitary $\cU:H\to H$ such that $\mA_1=\cU\mA_2\cU^{-1}$. If $\cA_1$ and $\cA_2$ are unitary equivalent then they belong to the same equivalence class from $\mathfrak{F}$ defined in Definition \ref{D2}.
\end{theorem}

While Theorem \ref{Fiso} is trivial we leave it for possible future improvements. One of the improvement is to rewrite the equivalence condition in $\mathfrak{F}$ in such a way that the inverse statement also becomes true. Namely: when is there $1-1$ mapping between the classes $\mathfrak{F}$ and unitary equivalent AF algebras. This statement is true for finite-dimensional $C^*$-algebras because Hilbert spaces of the same dimension are unitary isomorphic and we have finite partition of the Hilbert space onto orthogonal sum of its Hilbert subspaces. For infinite-dimensional AF-algebras there are counterexamples.

{\bf Example.} Let $H$ be a separable Hilbert space. Let $\mK$ be $C^*$-algebra of compact operators acting on $H$. Consider an one-dimensional extension $H_1=\C\os H$ and $C^*$-algebra $\wt\mK=0\os \mK$. It is seen that the corresponding unitalizations $\wt\mK_1$ and $\mK_1$ belong to the same class of equivalence in $\mathfrak{F}$ but they are not unitary equivalent. This is because for any $e\in H$ there is non-compact $\cK\in\mK_1$ such that $\cK e=0$. The last statement is not true for $\wt\mK_1$ since vectors $e\ne0$ from the one-dimensional supplement to $H$ are not null vectors for non-compact operators from $\wt\mK_1$.


{\section{Main results}\lb{sec3}}

Let $N,M\in\N$ be positive integers. Let $L^2_{N,M}=L^2(\T^N\to\C^M)$ be the Hilbert space of periodic vector valued functions defined on the multidimensional torus $\T^N$, where $\T=\R/\Z\simeq[0,1)$. Everywhere in the article, it is assumed the Lebesgue measure in the definition of Hilbert spaces of square-integrable functions. Let $R^{\iy}_{N,M}=R^{\iy}(\T^N\to\C^{M\ts M})$ be the $C^*$-algebra of matrix-valued regulated functions with rational discontinuities. The regulated functions with possible rational discontinuities are the functions that can be uniformly approximated by the step functions of the form
\[\lb{REG}
 {\bf S}({\bf x})=\sum_{n=1}^P \chi_{J_n}({\bf x}){\bf S}_{n},
\]
where $P\in\N$, ${\bf S}_{n}\in\C^{M\ts M}$, and $\chi_{J_n}$ is the characteristic function of the parallelepiped $J_n=\prod_{i=1}^N[p_{in},q_{in})$ with rational end points $p_{in},q_{in}\in\mathbb{Q}/\Z\ss\T$. In particular, continuous matrix-valued functions belong to $R^{\iy}_{N,M}$. Let us introduce the generating operators for integro-differential algebras. These operators are operators of multiplication by a function $\cM$, finite differential operators $\cD$, and integral operators $\cI$, all of them act on $L^2_{N,M}$: 
\[\lb{bas}
 \begin{array}{l} \cM_{\bf S}{\bf u}({\bf x})={\bf S}({\bf x}){\bf u}({\bf x}) \\ \cD_{i,h}{\bf u}({\bf x})=h^{-1}({\bf u}({\bf x}+h{\bf e}_i)-{\bf u}({\bf x}))\\ \cI_i{\bf u}=\int_0^1{\bf u}({\bf x})dx_i \end{array},\ \ \ \ {\bf u}({\bf x})\in L^2_{N,M},\ \ {\bf x}\in\T^N,
\] 
where the function ${\bf S}\in R^{\iy}_{N,M}$, the index $i\in\N_N=\{1,...,N\}$, the step of differentiation $h\in\mathbb{Q}$, the standard basis vector ${\bf e}_i=(\d_{ij})_{j=1}^N$, and $\d_{ij}$ is the Kronecker symbol. The $C^*$-algebra of finite-integro-differential operators is generated by all the operators \er{bas}
\[\lb{alg}
 \mF_{N,M}=\ol{\rm Alg}^{\mB}\{\cM_{\bf S},\ \cD_{i,h},\ \cI_i:\ \ {\bf S}\in R^{\iy}_{N,M},\ i\in\N_N,\ h\in\mathbb{Q}\},
\]
where $\mB\ev\mB_{N,M}=\mB(L^2_{N,M})$ is the $C^*$-algebra of all the bounded operators acting on $L^2_{N,M}$. The typical example of an operator $\cA$ from $\mF_{1,1}$ is
$$
 \cA u(x)=\sum_{n=1}^pA_n(x)\cD_{1,\frac np} u(x)+\int_0^1K(x,y)u(y)dy,\ \ u\in L^2_{1,1},\ x\in\T,
$$ 
where $A_n\in R^{\iy}_{1,1}$, $K\in R^{\iy}_{2,1}$, and $p\in\N$.
Let us provide the characterization of $\mF_{N,M}$.

\begin{theorem}\lb{T3}
The AF-algebra $\mF_{N,M}$ has the following Glimm-Bratteli symbol
\[\lb{H11}
 \mathfrak{n}(\mF_{N,M})=\lt(\prod_{n=2}^{\iy}\ma n & 0 \\ n-1 & 1 \am^{\otimes N}\rt)\ma 1 \\ 1 \am^{\otimes N}.
\]
In particular, $\mF_{N,M}$ and $\mF_{N_1,M_1}$ are isomorphic if and only if $N=N_1$.
\end{theorem}

Integro-differential algebras with different number of variables are non-isomorphic. This fact distinguishes these algebras from the differential algebras $\mH_{N,M}$ generated by $\cM_{\bf S}$ and $\cD_{i,h}$. The algebras $\mH_{N,M}$ are isomorphic to the universal UHF-algebra $\mU=\bigotimes_{n=1}^{\iy}\C^{n\ts n}$ independently on the number of variables $N$ and the number of functions $M$, see \cite{K2018arxiv}.

%

{\bf Example.} Let us consider the $C^*$-algebra of two-dimensional integro-differential operators $\mF_{2,M}$. We have
$$
 \prod_{i=1}^n\ma i & 0 \\ i-1 & 1 \am^{\otimes 2}\ma 1 \\ 1\am^{\otimes 2}=\ma n! \\ n! \am^{\otimes2}=(n!)^2\ma 1 \\ 1 \\ 1 \\ 1 \am
$$
and
$$
 \ma n+1 & 0 \\ n & 1 \am^{\otimes 2}=\ma (n+1)^2 & 0 & 0 & 0 \\ n(n+1) & n+1 & 0 & 0 \\
   n(n+1) & 0 & n+1 & 0 \\
   n^2 & n & n & 1 \am.
$$
Hence, the fragment of the Bratteli diagram for $\mF_{2,M}$ is

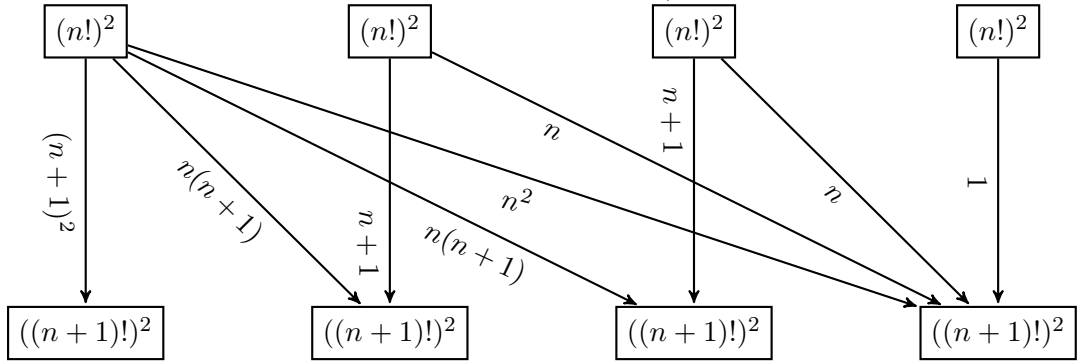
\begin{tikzpicture}[->,>=stealth',shorten >=1pt,auto,node distance=4cm,
                    thick,main node/.style={rectangle,draw,font=\sffamily\small\bfseries}]

  \node[main node] (1) {$(n!)^2$};
  \node[main node] (2) [right of=1] {$(n!)^2$};
  \node[main node] (3) [right of=2] {$(n!)^2$};
  \node[main node] (4) [right of=3] {$(n!)^2$};
  \node[main node] (5) [below of=1] {$((n+1)!)^2$};
  \node[main node] (6) [right of=5] {$((n+1)!)^2$};
  \node[main node] (7) [right of=6] {$((n+1)!)^2$};
  \node[main node] (8) [right of=7] {$((n+1)!)^2$};
  
  \path[every node/.style={font=\sffamily\small}]
   
   (1) edge node [below,rotate=270] {$(n+1)^2$} (5)
       edge node [below, rotate=315] {$n(n+1)$} (6)
       edge node [pos=0.7, below, rotate=335] {$n(n+1)$} (7)
       edge node [below, rotate=345] {$n^2$}    (8)
   (2) edge node [pos=0.75,below,rotate=270] {$n+1$} (6)
       edge node [pos=0.25, below,rotate=335] {$n$} (8)
   (3) edge node [pos=0.25,below,rotate=270] {$n+1$} (7)
       edge node [below, rotate=315] {$n$} (8)
   (4) edge node [below,rotate=270] {$1$} (8);
  
\end{tikzpicture}

Here, the vertices in the row represent direct summands of the finite dimensional sub-algebra, the edges represent partial embeddings into the next finite-dimensional sub-algebra appearing in the direct limit, and the edge labels are multiplicities of partial embeddings.

{\bf Remark.} Let us consider the algebra of one-dimensional scalar integro-differential operators $\mF_{1,1}$. The E-matrices for $\mF_{1,1}$ are given by Theorem \ref{T3}
$$
 {\bf E}_0=\ma 1\\ 1\am,\ \ {\bf E}_n=\ma n+1 & 0 \\ n & 1 \am.
$$
It is clear that
$$
 {\bf E}_n=\ma n+1 & 0 \\ 0 & 1 \am\ma 1 & 0 \\ 1 & 1 \am^n.
$$
Thus, there are arbitrary large sequences $\ma 1 & 0 \\ 1 & 1 \am^n$, $n\in\N$ in the direct limit for $\mF_{1,1}$. Remembering that these sequences correspond to the unitalized algebra of compact operators $\C1+\mK(L^2_{1,1})$, see above and, e.g., Example 3.3.1 in \cite{N2016}, we can expect that  $\mK(L^2_{1,1})\ss\mF_{1,1}$. This is true because any compact operator can be uniformly approximated by finite-dimensional operators in some orthonormal basis of $L^2_{1,1}$. Taking Walsh basis $f_n$, $n\in\N$ consisting of step functions, we see that for any $n,m\in\N$ the one-rank operator $\cC_{n,m}$ given by $u\to f_m\int_0^1f_nu$, where $u\in L^2_{1,1}$, belongs to $\mF_{1,1}$. Hence, any compact operator belongs to $\mH_{1,1}$, since it can be uniformly approximated by linear combinations of $\cC_{n,m}$. 

Finally, note that $\ma n+1 & 0 \\ 0 & 1 \am=(n+1)\os(1)$. E-matrices $(n+1)$, $n\in\N$ correspond to the universal uniformly hyper-finite algebra $\mU=\bigotimes_{n=1}^{\iy}\C^{n\ts n}$ which has the supernatural number $\mathfrak{n}(\mU)=\prod_{n=1}^{\iy}n$. Generated by $\cM_{\bf S}$ and $\cD_{i,h}$, see \er{bas}, $\mU$ is a sub-algebra of $\mF_{1,1}$. Roughly speaking, $\mF_{1,1}$ is a combination of  the universal UHF-algebra $\mU$ and the algebra of compact operators $\mK$. 

The natural extension of $\mF_{1,1}$ (or $\mF_{1,M}$) is the AF-algebra $\mF_1$ generated by the following commutative semigroup
$$
 \mathfrak{n}(\mF_{1})=\prod_{n=1}^{\iy}\prod_{m=1}^n\ma n & 0 \\ n-m & m \am.
$$
This is the maximal commutative semigroup of $2\ts2$-matrices from $\mathfrak{E}$ having the eigenvectors $\ma 1 \\ 1 \am$ and $\ma 0 \\ 1 \am$. Perhaps, it would be interesting to see the "physical meaning" of extended integro-differential operators from $\mF_1$.


{\section{Proof of the main results}\lb{sec4}}

{\bf Proof of Theorem \ref{T2}} The conditions of Definition \ref{D1} will be checked. We set ${\bf C}_0={\bf A}_0$, ${\bf C}_1={\bf B}_1$, and $r_1=1$, $m_1=2$ correspondingly. Next, ${\bf B}_1={\bf A}_{n_1}$ for some $n_1\ge1$. We take ${\bf C}_2=\prod_{i=0}^{n_1-1}{\bf A}_i$, or ${\bf C}_2={\bf A}_2$ if $n_1=1$. In the first case we set $r_2=n_1+1$, in the second case we set $r_2=3$.

Anyway, ${\bf C}_2\in\{{\bf A}_n\}_{n=1}^{\iy}$, since this is the semigroup. Hence, for some ${\bf C}_3\in\{{\bf A}_n\}_{n=1}^{\iy}$, we have $({\bf B}_1{\bf C}_2){\bf C}_3=\prod_{i=1}^{m_2-1}{\bf B}_i$ by the condition ($\pmb{\s}$). Thus ${\bf C}_2{\bf C}_3=\prod_{i=m_1}^{m_2-1}{\bf B}_i$ because ${\bf B}_1$ is invertible and all the matrices are commute.

By induction, suppose that for some $n>1$ we already found $1=r_1<...<r_n$, and $2=m_1...<m_n$, and ${\bf C}_i\in\{{\bf A}_n\}_{n=1}^{\iy}$ satisfying
\[\lb{i1}
 {\bf C}_{2j-1}{\bf C}_{2j-2}=\prod_{i=m_{j-1}}^{m_{j}-1}{\bf B}_i,\ \ {\bf C}_{2j-2}{\bf C}_{2j-3}=\prod_{i=r_{j-1}}^{r_{j}-1}{\bf A}_i,\ \ 2\le j\le n.
\]
Let $\mu({\bf A})$ be the maximal element of the matrix ${\bf A}$. It is true 
\[\lb{ineqmu}
 \mu({\bf A}_n{\bf A}_m)\ge\max(\mu({\bf A}_n),\mu({\bf A}_m)),
\] 
since ${\bf A}_n$, ${\bf A}_m$ are matrices with non-negative integer entries, without zero rows and columns. There are two possibilities: (a) $\lim_{p\to\iy}\mu({\bf C}_{2n-1}^p)=\iy$, and (b) $\mu({\bf C}_{2n-1}^p)$ are bounded. In the case (a), for some sufficiently large $p>1$ we have ${\bf C}_{2n-1}^p={\bf A}_r$, where $r>r_n$. We set $r_{n+1}=r+1$, ${\bf C}_{2n}={\bf C}_{2n-1}^{p-1}\prod_{i=r_n}^{r-1}{\bf A}_i$.
Hence, we obtain
\[\lb{r1}
 {\bf C}_{2n}{\bf C}_{2n-1}=\prod_{i=r_n}^{r_{n+1}-1}{\bf A}_i.
\]
Note that ${\bf C}_{2n}\in\{{\bf A}_n\}_{n=1}^{\iy}$, since this is the semigroup. Another possibility: (b) $\mu({\bf C}_{2n-1}^p)$ are uniformly bounded for all $p$. Then ${\bf C}_{2n-1}^p={\bf C}_{2n-1}^s$ for some $p> s$ because $\{{\bf C}_{2n-1}^p\}$ is a sequence of matrices with bounded non-negative integer entries. The existence of inverse matrix ${\bf C}_{2n-1}^{-1}$ leads to ${\bf C}_{2n-1}^{p-s}={\bf I}$ is the identity matrix. We set $r_{n+1}=r_n+1$, ${\bf C}_{2n}={\bf C}_{2n-1}^{p-s-1}{\bf A}_{r_n}$. These values also satisfy \er{r1}. Note that ${\bf E}$-matrices satisfying the condition (b) correspond to a permutation of elements in the Bratteli diagrams.

Again, there are two possibilities: (a) $\lim_{p\to\iy}\mu({\bf C}_{2n-1}^p)=\iy$, and (b) $\mu({\bf C}_{2n-1}^p)$ are bounded. Consider the first case (a), the second (b) can be treated as above. There is $p\ge1$ such that
\[\lb{clarge}
 \mu({\bf C}_{2n}^p)>\mu(\prod_{i=1}^{m_n}{\bf B}_i).
\]
Hence, by the condition ($\pmb{\s}$), taking ${\bf A}_p=(\prod_{i=1}^{m_n-1}{\bf B}_i){\bf C}_{2n}^p$ (recall that the set $\{{\bf A}_n\}_{n=1}^{\iy}$ is a semigroup) we have
\[\lb{mn}
 ((\prod_{i=1}^{m_n-1}{\bf B}_i){\bf C}_{2n}^p){\bf A}_r=\prod_{i=1}^{m_{n+1}-1}{\bf B}_i
\]
for some $m_{n+1}>m_n$ because of \er{clarge} and \er{ineqmu}.
We set ${\bf C}_{2n+1}={\bf C}_{2n}^{p-1}{\bf A}_r$. Using \er{mn}, we deduce
\[\lb{m1}
 {\bf C}_{2n+1}{\bf C}_{2n}=\prod_{i=m_n}^{m_{n+1}-1}{\bf B}_i.
\]
Thus, by induction we prove that $\prod_{n=0}^{\iy}{\bf A}_n$ and $(\prod_{n=1}^{\iy}{\bf B}_n){\bf A}_0$ are equivalent, see Definition \ref{D1}. By Theorem \ref{Tstruc}, they represent the same algebra.


If $\mu(\prod_{i=1}^p{\bf B}_i)$ are bounded for all $p$ then 
there is ${\bf A}_r$ and $1\le m_1<m_2<...$ such that
$
 \prod_{i=1}^{m_n}{\bf B}_i={\bf A}_r
$
for all $n$. Thus, $\mathfrak{n}^{-1}((\prod_{n=1}^{\iy}{\bf B}_n){\bf A}_0)\cong\mM({\bf A}_r{\bf A}_0)$ is a sub-algebra of $\mM(\prod_{n=0}^r{\bf A}_n)$, which, in turn, the sub-algebra of $\mathfrak{n}^{-1}(\prod_{n=0}^{\iy}{\bf A}_n)$. Now, suppose that $\mu(\prod_{i=1}^p{\bf B}_i)\to\iy$. Then we can take $1=m_1<m_2<...$ such that
$$
 \prod_{i=m_n}^{m_{n+1}-1}{\bf B}_i={\bf A}_{r_n}\ \ n\ge1,
$$
where $0=r_0<r_1<r_2<....$ Denote 
$$
 {\bf D}_n=\prod_{i=r_{n-1}+1}^{r_n}{\bf A}_i,\ \ {\bf E}_n=\prod_{j=1}^n(\prod_{i=r_{j-1}+1}^{r_j-1}{\bf A}_i),
$$
where ${\bf E}_n={\bf I}$ is the identity matrix if $r_n-2<r_{n-1}$. Then the following infinite commutative diagrams
$$
  \begin{tikzcd}
    \mM({\bf A}_0) \arrow{r}{{\bf D}_1}  & \mM({\bf D}_1{\bf A}_0) \arrow{r}{{\bf D}_2}  & \mM({\bf D}_2{\bf D}_1{\bf A}_0) & ...\arrow{r}{} & \mathfrak{n}^{-1}(\prod_{n=0}^{\iy}{\bf A}_n) \\
    \mM({\bf A}_0)\arrow{r}{{\bf A}_{r_1}}\arrow{u}{{\bf I}} & \mM({\bf A}_{r_1}{\bf A}_0)\arrow{u}{{\bf E}_1}\arrow{r}{{\bf A}_{r_2}} & \mM({\bf A}_{r_2}{\bf A}_{r_1}{\bf A}_0)\arrow{u}{{\bf E}_2} & ...\arrow{r}{} & \mathfrak{n}^{-1}((\prod_{n=1}^{\iy}{\bf B}_n){\bf A}_0)\arrow{u}{} 
  \end{tikzcd}
$$
show that $\mathfrak{n}^{-1}((\prod_{n=1}^{\iy}{\bf B}_n){\bf A}_0)$ is the sub-algebra of $\mathfrak{n}^{-1}(\prod_{n=0}^{\iy}{\bf A}_n)$. \BBox

{\bf Proof of Theorem \ref{Tiso}.} Suppose that $N>M$. If $\mathfrak{n}^{-1}(\prod_{n=0}^{\iy}{\bf A}_n)\cong\mathfrak{n}^{-1}(\prod_{n=0}^{\iy}{\bf B}_n)$ then there is the sequence of matrices $\{{\bf C}\}_{n=0}^{\iy}$ satisfying \er{D11}, namely
$$
 {\bf C}_{2n-2}{\bf C}_{2n-3}=\prod_{i=r_{n-1}}^{r_{n}-1}{\bf A}_i,\ \ {\bf C}_{2n-1}{\bf C}_{2n-2}=\prod_{i=m_{n-1}}^{m_{n}-1}{\bf B}_i,\ \ {\bf C}_{2n}{\bf C}_{2n-1}=\prod_{i=r_{n}}^{r_{n+1}-1}{\bf A}_i
$$
for some $n>2$. This yields to
$$
 \prod_{i=r_{n-1}}^{r_{n+1}-1}{\bf A}_i={\bf C}_{2n}(\prod_{i=m_{n-1}}^{m_{n}-1}{\bf B}_i){\bf C}_{2n-3}. 
$$
The matrix in LHS has the full rank $N$, while the matrix in RHS has a rank not more than $M$. This is the contradiction.
\BBox

{\bf Proof of Theorem \ref{T3}.} Let us start from the 1D case $N=M=1$. For $h\in\mathbb{Q}$, define the shift operator $\cS_{h}=1-h\cD_{1,h}$. Define also the operators of multiplication by the characteristic functions of intervals
\[\lb{Mip}
 \cM_{j,p}\ev\cM_{\chi_{I_j^p}},\ \ \ I_j^p=\lt[\frac{j}p,\frac{j+1}p\rt),\ \ \ j\in\Z_p=\{0,...,p-1\},\ \ \ p\in\N.
\]
The operators satisfy some elementary properties
\begin{multline}\lb{prop}
 \cM_{i,p}\cM_{j,p}=\d_{ij}\cM_{i,p},\ \ \ \cS_{\frac jp}\cM_{i,p}=\cM_{i+j,p}\cS_{\frac jp},\ \ \ \cS_h\cS_t=\cS_{h+t},\ \ \ \cM_{i,p}^*=\cM_{i,p},\\ 
 \cS_{h}^*=\cS_{-h},\ \ \ \cS_{h}\cI_1=\cI_1\cS_h=\cI_1,\ \ \ \cI_1\cM_{i,p}\cI_1=p^{-1},
\end{multline}
where $i,j\in\Z_p$, $h,t\in\mathbb{Q}$, and $p\in\N$. For $i,j\in\Z_p$, define the basis operators
\[\lb{basop}
 \cB_{i,j}^p=p\cM_{i,p}\cI_1\cM_{j,p},\ \ \ \cA_{i,j}^p=\cM_{i,p}\cS_{\frac{i-j}p}-\cB_{i,j}^p.
\]
Using \er{prop}, we can directly check the properties
\[\lb{matprop}
 \cB_{i,j}^p\cB_{n,m}^p=\d_{jn}\cB^p_{i,m},\ \ \ (\cB^p_{i,j})^*=\cB^p_{j,i},\ \ \ 
 \cA^p_{i,j}\cA^p_{n,m}=\d_{jn}\cA^p_{i,m},\ \ \ (\cA_{i,j}^p)^*=\cA_{j,i},\ \ \ \cA_{i,j}\cB_{n,m}=0.
\]
Identities \er{matprop} means that
\[\lb{matalg}
 \mH_p\ev{\rm Alg}\{\cA_{i,j}^p,\ \cB_{i,j}^p:\ \ i,j\in\Z_p\}\cong\mM(p)\os\mM(p)=\mM(p,p)
\]
with the $*$-isomorphism defined by
\[\lb{matiso}
 \cA_{ij}^p\mapsto(\d_{in}\d_{jm})_{n,m=0}^{p-1}\os{\bf 0}_p,\ \ \ \cB_{ij}^p\mapsto{\bf 0}_p\os(\d_{in}\d_{jm})_{n,m=0}^{p-1},
\]
where ${\bf 0}_p$ is the zero element in $\mM_p$. Let $q\in\N$ be some positive integer. Using \er{basop} and the identity
$$
 \cM_{i,p}=\sum_{n=iq}^{(i+1)q-1}\cM_{n,pq},
$$
we obtain
\begin{multline}\lb{Bemb}
 \cB_{ij}^p=\frac1q\sum_{n=iq}^{(i+1)q-1}\sum_{m=jq}^{(j+1)q-1}\cB_{nm}^{pq},\ \ \ \cA_{ij}^p=(\sum_{n=iq}^{(i+1)q-1}\cM_{n,qp})\cS_{\frac{iq-pq-n+n}{pq}}-\cB_{i,j}^p=\\
 \sum_{n=iq}^{(i+1)q-1}(\cA^{pq}_{n+(j-i)q}+\cB^{pq}_{n+(j-i)q})-\frac1q\sum_{n=iq}^{(i+1)q-1}\sum_{m=jq}^{(j+1)q-1}\cB_{nm}^{pq}.
\end{multline}
Identity \er{Bemb} shows how $\mH_{p}$ is embedded into $\mH_{pq}$. Namely, the corresponding $*$-embedding is defined by
\[\lb{Embq}
 {\bf A}\os{\bf B}\mapsto({\bf A}\otimes{\bf I}_q)\os({\bf A}\otimes{\bf I}_q-{\bf A}\otimes(q^{-1}{\bf 1}_q)+{\bf B}\otimes q^{-1}{\bf 1}_q),
\]
where ${\bf I}_q$ is the identity matrix in $\mM_q$ and ${\bf 1}_q=(1)\in\mM_q$ is the matrix, which has all entries equal to $1$. The matrix $q^{-1}{\bf 1}_q$ is the rank-one matrix and, hence, it is unitarily equivalent to the matrix with one non-zero entry
\[\lb{1q}
 q^{-1}{\bf 1}_q\simeq\ma 1 & 0 & ... \\
                          0 & 0 & ... \\
                          ... & ... & ... \am.                          
\] 
Using \er{1q} we conclude that $*$-embedding \er{Embq} between $\mH_{p}\cong\mM(p,p)$ and $\mH_{pq}\cong\mM(pq,pq)$ has E-matrix
\[\lb{Emat1}
 {\bf E}=\ma q & 0 \\ q-1 & 1 \am.
\]
The integral operator $\cI_1$ belongs to all $\mH_p$, since
\[\lb{int1}
 \cI_1=\sum_{i\in\Z_p}\sum_{j\in\Z_p}\cB_{ij}^p\in\mH_p.
\]
By definition, any $S\in R^{\iy}_{1,1}$ can be uniformly approximated by step functions with rational discontinuities. Thus, the operator of multiplication by the function $\cM_{S}$ can be uniformly approximated by linear combinations of $\cM_{i,p}$, $i\in\Z_p$. On the other hand, using \er{basop}, we have
\[\lb{Mip1}
 \cM_{i,p}=\cA_{i,i}^p+\cB_{i,i}^p\in\mH_{p}
\]
Hence, for any $S\in R^{\iy}_{1,1}$, the operator $\cM_{S}$ can be uniformly approximated by the elements from $\mH_p$ with arbitrary precision when $p\to\iy$.
The identity operator $1$ belongs to all $\mH_p$, since
\[\lb{id1}
 1=\sum_{j=0}^{p-1}\cM_{j,p}.
\]
The shift operators $\cS_h$ (with $h=q/p\in\mathbb{Q}$) belongs to $\mH_p$, since
\[\lb{shift1}
 \cS_{\frac1p}=\sum_{i\in\Z_p}(\cA_{i,i-1}^p+\cB_{i,i-1}^p)\in\mH_p,\ \ \cS_{\frac qp}=\cS_{\frac 1p}^q\in\mH_p
\]
by \er{basop}, \er{prop}, and \er{id1}. Hence, the finite differentials belong also to $\mH_p$: 
\[\lb{diff1}
 \cD_{1,h}=h^{-1}(1-\cS_{h})\in\mH_p.
\]
Using \er{int1}, \er{diff1}, and the mentioned above fact about the approximation of $\cM_S$ (for any $S\in R^{\iy}_{1,1}$) by the elements from $\mH_p$, we conclude that $\mF_{1,1}$ is the inductive limit of $\mH_p$ for $p\to\iy$. In particular, taking $p_n=n!$ and using \er{Emat1} for $\mH_{p_n}\ss\mH_{p_{n+1}}$ we obtain \er{H11} for $N=M=1$. 

The algebra $\mF_{1,M}=\mM(M)\otimes\mF_{1,1}$ has the same Glimm-Bratteli symbol as $\mF_{1,1}$, since
\begin{multline}
 \mathfrak{n}(\mF_{1,M})=\mathfrak{n}(\mM(M)\otimes\mF_{1,1})=\lt(\prod_{n=2}^{\iy}\ma n & 0 \\ n-1 & 1 \am\rt)\ma 1 \\ 1 \am M=\lt(\prod_{n=2}^{\iy}\ma n & 0 \\ n-1 & 1 \am\rt)\ma M \\ M \am=
\\
 \lb{FM1b}
 \lt(\prod_{n=2}^{\iy}\ma n & 0 \\ n-1 & 1 \am\rt)\ma M & 0 \\ M-1 & 1\am\ma 1 \\ 1 \am=\lt(\prod_{n=2}^{\iy}\ma n & 0 \\ n-1 & 1 \am\rt)\ma 1 \\ 1 \am=\mathfrak{n}(\mF_{1,1}).
\end{multline}
Let us discuss why the first identity in the last string of \er{FM1b} is true. The matrices \er{Emat1} form a commutative (multiplicative) semigroup. Then, the infinite product with one duplicated term in LHS of \er{FM1b} obviously satisfies the condition ($\pmb{\s}$) from Theorem \ref{T2}. Hence, $\mF_{1,M}$ and $\mF_{1,1}$ are isomorphic. There is also a more intuitive similarity with supernatural numbers
$$
 \lt(\prod_{n=2}^{\iy}\ma n & 0 \\ n-1 & 1 \am\rt)\ma M & 0 \\ M-1 & 1\am\ma 1 \\ 1 \am=\lt(\prod_{n=1}^{\iy}\ma p_n & 0 \\ p_n-1 & 1 \am^{\iy}\rt)\lt(\prod_{j=1}^{K}\ma p_{n_j} & 0 \\ p_{n_j}-1 & 1 \am^{R_j}\rt)\ma 1 \\ 1 \am=
$$
$$
 \prod_{n=1}^{\iy}\ma p_n & 0 \\ p_n-1 & 1 \am^{\iy}\ma 1 \\ 1 \am=\prod_{n=2}^{\iy}\ma n & 0 \\ n-1 & 1 \am\ma 1 \\ 1 \am,
$$
where $p_1=2$, $p_2=3$, $p_3=5$, ... are prime numbers and $M=\prod_{j=1}^K p_{n_j}^{R_j}$ is the prime factorization of $M$. This similarity with supernatural numbers is possible because all the matrices are commute.

Consider the case $N>1$. Using the fact that $L^2_{N,M}=\bigoplus_{j=1}^M(L^2_{1,1})^{\otimes N}$ we deduce that $\mF_{N,M}=\mM(M)\otimes\mF_{1,1}^{\otimes N}$. This means that
$$
 \mathfrak{n}(\mM(M)\otimes\mF_{1,1}^{\otimes N})=\mathfrak{n}(\mF_{1,M}\otimes\mF_{1,1}^{\otimes N-1})=\mathfrak{n}(\mF_{1,1}^{\otimes N})=\lt(\prod_{n=2}^{\iy}\ma n & 0 \\ n-1 & 1 \am^{\otimes N}\rt)\ma 1 \\ 1 \am^{\otimes N}
$$
which proves \er{H11}. Thus, the $C^*$-algebras $\mF_{N_1,M_1}$ and $\mF_{N_2,M_2}$ are isomorphic if and only if $N_1=N_2$ by Theorem \ref{Tiso}. \BBox

\medskip

\bibliography{bibl_perp1}

\end{document}